\documentstyle[amstex,amssymb]{article}
\long\def\comment#1{}
\newcommand\calH{{\cal H}}
\newcommand\antisym{\mathop{\rm Ant}\nolimits}
\newcommand\sym{\mathop{\rm Sym}\nolimits}
\renewcommand\Im{\mathop{\rm Im}\nolimits}
\newcommand\End{\mathop{\rm End}\nolimits}
\newcommand\Id{\mathop{\rm Id}\nolimits}
\newcommand\isom{\cong}
\newcommand\tensor{\otimes}
\newcommand\Complex{{\Bbb C}}

\newcommand\Integer{{\Bbb Z}}
\newcommand\sh{\mathop{\mathrm sinh}\nolimits}
\newcommand\ch{\mathop{\mathrm cosh}\nolimits}
\newcommand\tr{\mathop{\mathrm tr}\nolimits}
\newcommand\qdet{\text{q-det\,}}
\newcommand\e{{\mathrm e}}
\renewcommand\i{{\mathrm i}}
\newcommand{\secref}[1]{\S\ref{#1}}

\begin{document}
\begin{titlepage}

\title{Bethe Ansatz for Higher Spin XYZ Models\\
       --- Low-lying Excitations ---}

\author{Takashi Takebe
\thanks{Present address: Department of Mathematics, the University of
California, Berkeley, CA94720, USA. (till August 1997).}\\
Department of Mathematical Sciences,\\
the University of Tokyo, \\
Komaba 3-8-1, Meguro-ku, \\
Tokyo, 153, Japan}

\maketitle

\abstract{
A higher spin generalization of the XYZ spin chain defined by the fusion
procedure is considered. Energy and momentum of the ground state and
low-lying excited states of the model are computed by means of an
isomorphism between a space of theta functions on which the Sklyanin
algebra acts and a space obtained by the fusion procedure.
}

\end{titlepage}
\section{Introduction.}
In this letter we calculate energy and momentum of low-lying excited
states of the higher spin generalization of the XYZ spin chain model
by means of the algebraic Bethe Ansatz method. This work is a
continuation of the author's previous works \cite{take:92},
\cite{take:95}, where a generalization of the eight-vertex model
\cite{bax} was studied and the two particle $S$-matrix of
the corresponding spin chain model was computed. There we assigned a
$(2l+1)$-dimensional space to the vertical edges of the lattice and a
two-dimensional space to the horizontal edges. In the present paper we
make use of the fusion procedure (see \cite{kul-resh-skl}, \cite{cher},
\cite{djkmo}, \cite{zhou-hou}) to assign a higher dimensional space to the
horizontal edges and then consider the corresponding one-dimensional
quantum spin chain model.

We shall give in \secref{fusion} an explicit isomorphism of a
representation space of the Sklyanin algebra which is defined as a space
of theta functions with the space of symmetric tensors. (This is a special
case of more general isomorphisms by Hasegawa \cite{has}.) Through this
identification we can identify the $L$ operator of the model in
\cite{take:95} with a special case of $R$ matrices in \cite{djkmo}, and
use the results in \cite{take:95} to compute the general models. In
particular this identification is indispensable to interpret a special
value of the logarithm of a transfer matrix and its logarithmic derivative
as a momentum and an energy operator. See, e.g.,
\cite{takh-fad:79}.

As is conjectured in \cite{take:95} from the corresponding results for the
higher spin XXX and XXZ models \cite{takh}, \cite{sog},
\cite{kir-resh}, the energy and the momentum are independent of the spin
of the local quantum space and expressed as a sum of two terms each of
which depends on a rapidity of a hole in the string configuration. Thus we
can justify the interpretation of those states as ``two-particle states''.

We make use of the results in \cite{take:95}, but adopt simpler
normalizations in \cite{take:96}. In particular, we normalize the $R$
matrices by the unitarity condition, and therefore they are meromorphic
functions of the spectral parameter, while \cite{djkmo} and
\cite{kir-resh} use holomorphic $R$ matrices.

\section{Fusion procedure.}
\label{fusion}
\setcounter{equation}{0}

In this section, we briefly review the fusion procedure for the elliptic
$R$ matrices.

Let $V^l = \sym(V_1 \tensor \cdots \tensor V_{2l})$ for
$l \in \frac{1}{2} \Integer$, where
$V_i \isom \Complex^2$ $(i = 1, \ldots, 2l)$ and $\sym$ is the
symmetrizer. The elliptic $R$ matrices $R^{l,l'}(u)$ have the following
properties. 

\begin{enumerate}
\renewcommand{\labelenumi}{(\roman{enumi})}
\item 
$R^{l,l'}(u)$ is a linear endomorphism of $V^l \tensor V^{l'}$
meromorphically depending on a complex parameter $u$.

\item({\em Yang-Baxter equation})
As an endomorphism of $V^l \tensor V^{l'} \tensor V^{l''}$,
\begin{multline}
    R_{12}^{l ,l' }(u_1 - u_2)
    R_{13}^{l ,l''}(u_1 - u_3)
    R_{23}^{l',l''}(u_2 - u_3) 
    =\\
    =
    R_{23}^{l',l''}(u_2 - u_3) 
    R_{13}^{l ,l''}(u_1 - u_3)
    R_{12}^{l ,l' }(u_1 - u_2).
\label{yb}
\end{multline}

\item({\em Unitarity})
As an endomorphism of $V^l \tensor V^{l'}$,
\begin{equation}
    R_{12}^{l,l'}(u - v) R_{21}^{l',l}(v - u) =
    \Id_{V^l\tensor V^{l'}}.
\label{unitarity}
\end{equation}

\item
When $u=0$, $R^{l,l}(u)$ is a permutation operator:
for all $v, w \in V^l$,
\begin{equation}
    R^{l,l}(0) (v \tensor w) = w \tensor v.
\label{R(0)=perm}
\end{equation}

\end{enumerate}

They are constructed in \cite{cher}, \cite{djkmo}, \cite{zhou-hou} by the
fusion procedure \cite{kul-resh-skl} from {\em Baxter's $R$ matrix}
$R^{1/2, 1/2}(u) = R(u;\tau)$ defined by
\begin{equation}
    R(u) = \sum_{a=0}^3 W_a(u) \sigma^a \tensor \sigma^a,\qquad
    W_a(u) := 
    \frac{\theta_{g_a}(u   ; \tau) \theta_{11}(    2\eta;\tau)}
       {2 \theta_{g_a}(\eta; \tau) \theta_{11}(u + 2\eta;\tau)},
\label{def:R}
\end{equation}
where $g_0 = (11)$, $g_1 = (10)$, $g_2 = (00)$, $g_3 = (01)$.
The explicit definition of $R^{l,l'}$ is:
\begin{align}
    R^{1/2,l'}_{V_i, V^{l'}}(u) 
    :=& \sym_{\overline{2l'} \ldots \overline{1}} 
    R_{V_i, V_{\overline{2l'}}}(u + (   2l'-1)\eta) \cdots
    \\ & \cdots
    R_{V_i, V_{\overline{\jmath}  }}(u + (2j-2l'-1)\eta) \cdots
    R_{V_i, V_{\overline{1     }  }}(u + (  -2l'+1)\eta),
\label{def:R-1/2-l'}
\\
    R^{l,l'}_{V^l, V^{l'}}(u)
    :=& \sym_{1 \ldots 2l}
    R^{1/2, l'}_{V_{2l}, V^{l'}} (u + (   2l-1)\eta) \cdots
    \\ &\cdots
    R^{1/2, l'}_{V_{ j}, V^{l'}} (u + (2j-2l-1)\eta) \cdots
    R^{1/2, l'}_{V_{ 1}, V^{l'}} (u + (  -2l+1)\eta).
\label{def:R-l-l'}
\end{align}
Here $V_i \isom V_{\overline{\jmath}} \isom \Complex^2$, suffixes of $R$
designate the spaces on which the $R$ matrix acts and $\sym_{1 \ldots m}$
is the symmetrizer on the space $V_1 \tensor \cdots \tensor V_m$ etc.

There is another expression of $R^{1/2,l'}(u)$ in terms of the
representation of the Sklyanin algebra which we used in
\cite{take:95}. The Sklyanin algebra \cite{skl:82}
$U_{\tau,\eta}(sl(2))$ is generated by four
generators $S^0$, $S^1$, $S^2$ and $S^3$ satisfying the relations coming
from the Yang-Baxter type relation of the $L$ operator, $L(u)$, defined by
\begin{equation}
    L(u) = \sum_{a=0}^3 W^L_a(u) \sigma^a \tensor S^a, \qquad
    W^L_a(u) = \frac{                     \theta_{g_a}(u)}
                    {2 \theta_{11}(2\eta) \theta_{g_a}(\eta)}.
\label{def:L}
\end{equation}
The Sklyanin algebra has a representation on a space of theta functions
\begin{gather}
    \rho^{(l)}: U_{\tau,\eta}(sl(2)) 
              \to \End_{\Complex}(\Theta^{4l+}_{00}),
\\
    \Theta^{4l+}_{00} =
    \left\{f(y)\text{: holomorphic on }\Complex \, \left| \,
     {{f(y+1) = f(-y) = f(y),}\atop
      {f(y+\tau)=\e^{-4l\pi \i(2y+\tau)}f(y)}}
    \right.\right\}.
\label{def:rep-space}
\end{gather}
It is easy to see that $\dim\Theta^{4l+}_{00} = 2l+1$. The generators
$S^a$ act on this space as difference operators. See \cite{skl:83} and
Appendix A. of \cite{take:95}, \cite{take:96}. We fix a basis of
$\Theta^{2+}_{00}$,
$(\theta_{00}(2y;2\tau)-\theta_{10}(2y;2\tau),
  \theta_{00}(2y;2\tau)+\theta_{10}(2y;2\tau))$
and identify $\Theta^{2+}_{00}$ with $\Complex^2$ through this basis.
Then we fix an isomorphism of the space of symmetric tensors and the spin
$l$ representation space of the Sklyanin algebra as follows:
\begin{multline}
    V^l = \sym( V_1 \tensor \cdots \tensor V_{2l}) \owns
    \sym( f_1 (y_1) \tensor \cdots \tensor f_{2l}(y_{2l}) )
    \mapsto \\
    \mapsto
    f_1 (y) \cdots f_{2l}(y) \in V^l = \Theta^{4l+}_{00},
\label{isom:sym-tensor=theta}
\end{multline}
where $V_i \isom \Complex^2$ $(i = 1, \ldots, 2l)$ identified with
$\Theta^{2+}_{00}$, $f_i(y_i) \in V_i$ and $\sym$ is the symmetrizer.
Under this identification the $L$ operator \eqref{def:L} is
proportional to the $R$ matrix defined by \eqref{def:R-1/2-l'}:
\begin{equation}
    R^{1/2,l} (u) 
    = \frac{\theta_{11}(2\eta)}{\theta_{11}(u+(2l+1)\eta)}
      \Id_{\Complex^2} \tensor \rho^{(l)}(L(u+\eta)).
\label{R=L}
\end{equation}
This can be verified by comparing the action of both hand sides on the
intertwining vectors. (See Lemma 2.1.3 and Theorem 2.3.3 of \cite{djkmo}
(Adv.\ Stud.\ Pure Math.\ {\bf 16}) for the left hand side and (1.18--21)
of \cite{take:96} for the right hand side. See also \cite{has} for general
cases.)

As is the case with the trigonometric and rational $R$ matrix, there is a
recurrence relation with respect to the auxiliary spin \cite{kul-resh},
\cite{kir-resh}. Let $R^{l,l'}(u)$ and $R^{1/2,l'}(u)$ are the $R$
matrices on the space $\sym(V_{2l} \tensor \cdots V_1) \tensor V^{l'}$ and
on $V_0 \tensor V^{l'}$ respectively. 
Here $V_i \isom \Complex^2$ $(i = 1, \ldots, 2l)$ and $V^{l'}$ is a space
of symmetric tensors defined above. Then as an operator on
$\sym(V_{2l} \tensor \cdots V_1) \tensor V_0 \tensor V^{l'}$,
\begin{multline}
    R^{l,l'}(u + \eta) R^{1/2,l'}(u - 2l\eta)
    = \\
    =
    \left(
    \begin{array}{c|c}
    R^{l+1/2,l'}(u) & 0 \\
    \hline 
    \ast          & \qdet R^{1/2,l'}(u-(2l-1)\eta) \times \\
                  & \hfill\times R^{l-1/2,l'}(u+2\eta)
    \end{array}
    \right),
\label{recurrence:R}
\end{multline}
where $\qdet R$ is the quantum determinant \cite{kul-skl} defined by
\begin{equation}
    \qdet R^{1/2,l'}(u) = 
    \tr_{01} P^-_{01} R^{1/2,l'}(u+\eta) R^{1/2,l'}(u-\eta)
    =
    \frac{\theta_{11}(u - 2 l' \eta)}{\theta_{11}(u + 2 l' \eta)}
    \Id_{V^{l'}},
\label{def:q-det}
\end{equation}
($P^-_{01}$ is a projection to the antisymmetric tensor in 
$V_1 \tensor V_0$). The block structure of the right hand side of
\eqref{recurrence:R} comes from the decomposition of the tensor product of
the auxiliary spaces by the Young symmetrizers as follows:
\begin{equation}
    \lower 1ex
    \vbox{\offinterlineskip
    \hrule
    \halign{&\strut\vrule\hfil#\hfil\cr
    \ 1 &\ 2 &\ $\cdots$ &\ $2l$ & \cr
    }
    \hrule}
    \tensor
    \lower 1ex
    \vbox{\offinterlineskip
    \hrule
    \halign{&\strut\vrule\hfil#\hfil\cr
    \ 0 &\cr
    }
    \hrule}
    =
    \lower 1ex
    \vbox{\offinterlineskip
    \hrule
    \halign{&\strut\vrule\hfil#\hfil\cr
    \ 0 &\ 1 &\ 2 &\ $\cdots$ &\ $2l$ & \cr
    }
    \hrule}
    \oplus
    \lower 3.9ex
    \vbox{\offinterlineskip
    \hrule
    \halign{&\strut\vrule\hfil#\hfil\cr
    \ 1 &\ 2 &\ $\cdots$ &\ $2l$ & \cr
    \noalign{\hrule}
    \ 0 & & \omit\cr
    \omit\hrulefill&\omit\cr
    }
    }
\label{decomposition}
\end{equation}
Here we denote the image of the Young symmetrizer by the corresponding
Young tableau. Since the second Young symmetrizer in the right hand side
of \eqref{decomposition} gives an isomorphism from the space 
$\sym(V_{2l} \tensor \cdots \tensor V_2) \tensor \antisym(V_1 \tensor V_0)
\isom \sym(V_{2l} \tensor \cdots \tensor V_2)$ to its image, we identify
these spaces in the right lower corner of the right hand side of
\eqref{recurrence:R}. The proof of this equation is based on Lemma 5 of
\cite{cher} and done by the same argument as that in section 1 of
\cite{kul-resh}.

\section{Higher spin XYZ model.}
\label{XYZ-model}
\setcounter{equation}{0}

In this section we define a higher spin generalization of the XYZ model
and apply the algebraic Bethe Ansatz.

The state space of our model is
\begin{equation}
    \calH = V^{l}_{\overline{1}} \tensor V^{l}_{\overline{2}}
            \tensor \cdots \tensor V^{l}_{\overline{N}},
\label{def:state-space}
\end{equation}
where $V^{l}_{\overline{\jmath}} \isom V^{l}$. We define the model by
specifying the transfer matrix, namely, generating function of quantum
integrals of motion, as follows:
\begin{equation}
    T^{l',l}(u) :=
    \tr_{V^{l'}} R^{l',l}_{V^{l'}, V^{l}_{\overline{N}}}(u) \cdots
                 R^{l',l}_{V^{l'}, V^{l}_{\overline{1}}}(u)
    \in \End_\Complex(\calH),
\label{def:T}
\end{equation}
where $R_{V^{l'}, V^l_{\overline\jmath}}$ acts non-trivially only on the
component $V^{l'} \tensor V^l_{\overline{\jmath}}$ of 
$V^{l'} \tensor \calH$. The most important property of the transfer
matrix is the commutativity
\begin{equation}
    [T^{l',l}(u), T^{l'',l}(v)] = 0,
\end{equation}
which is a consequence of the Yang-Baxter equation \eqref{yb}.

Thanks to \eqref{R(0)=perm}, we can define a momentum operator $p$ and a
Hamiltonian $H$ of the spin chain by:
\begin{equation}
    p = \frac{1}{\i} \log T^{l,l}(0), \qquad
    H = \text{const.} 
        \left. \frac{d}{du} \log T^{l,l}(u) \right|_{u=0},
\label{def:p,H}
\end{equation}
where const.\ is a constant which we do not fix here. (See
\cite{takh-fad:79}, \cite{takh-fad:81}, \cite{takh}, \cite{sog},
\cite{kir-resh}.) 

Due to the recurrence relation of the $R$ matrix \eqref{recurrence:R}, the
transfer matrix satisfies the following recurrence relation:
\begin{multline}
    T^{l'+1/2,l}(u) = \\
    =
    T^{l',l}(u+\eta) T^{1/2,l}(u-2l'\eta)
    - \frac{\theta_{11}(u+(-2l'+1-2l)\eta)}
           {\theta_{11}(u+(-2l'+1+2l)\eta)} T^{l'-1/2,l}(u+2\eta).
\label{recurrence:t}
\end{multline}
Hence the diagonalization problem of $T^{l',l}(u)$ reduces to that of
$T^{1/2,l}(u)$.

Hereafter we assume that the elliptic modulus $\tau$ is a pure imaginary
number $\tau = \i/t$, $t>0$ and the anisotropy parameter $\eta = r'/r$ is
a rational number. We also assume that $M = Nl$ is an integer.

It was shown in \cite{take:92} that there exist
vectors $\Psi_\nu(w_1, \ldots, w_M)$ depending on an integer $\nu$ and
complex parameters $(w_1, \ldots, w_M)$ such that
\begin{equation}
    T^{1/2, l}(u) \Psi_\nu(w_1, \ldots, w_M)
    =
    t^{1/2, l}(u; \nu, w_1, \ldots, w_M)
    \Psi_\nu(w_1, \ldots, w_M),
\label{bethe-vector}
\end{equation}
provided that $(\nu; w_1, \ldots, w_M)$ satisfy the Bethe equations:
\begin{equation}
    \left(\frac{\theta_{11}(w_j + 2l\eta;\tau)}
               {\theta_{11}(w_j - 2l\eta;\tau)}\right)^N
    =
    \e^{-4\pi \i \nu \eta}
    \prod\begin{Sb}k=1 \\ k \neq j \end{Sb}^M
    \frac{\theta_{11}(w_j - w_k + 2\eta;\tau)}
         {\theta_{11}(w_j - w_k - 2\eta;\tau)},
\label{bethe-eq}
\end{equation}
for all $j=1, \dots, M$. The eigenvalue $t^{1/2,l}(u)$ is equal to
\begin{equation}
    t^{1/2,l}(u;\nu, w_1, \ldots, w_M)
    =
    \frac{Q(u-2\eta)}{Q(u)} + h(u) \frac{Q(u+2\eta)}{Q(u)},
\label{bethe-ev}
\end{equation}
where
\begin{equation}
    Q(u):= \e^{-\pi\i \nu u} \prod_{j=1}^M \theta_{11}(u-w_j+\eta), 
    \qquad
    h(u):=\left(
          \frac{\theta_{11}(u+(-2l+1)\eta)}{\theta_{11}(u+(2l+1)\eta)}
          \right)^N.
\label{def:Q,h}
\end{equation}
Inductively using the recurrence relation \eqref{recurrence:t}, we can
prove that the eigenvalue of
$T^{l',l}(u)$ is
\begin{multline}
    t^{l',l}(u; \nu, w_1, \ldots, w_M) =
    \\=
    \sum_{j=0}^{2l'} a_j^{l',l}(u)
    \frac{Q(u+(2l'+1)\eta)    Q(u-(2l'+1)\eta)}
         {Q(u+(2l'+1-2j)\eta) Q(u+(2l'-1-2j)\eta)},
\label{bethe-ev:l'}
\end{multline}
where $a_0^{l',l}(u) = 1$,
$a_j^{l',l}(u) = \prod_{k=1}^j h(u+(2l'+1-2k)\eta)$.

\section{Thermodynamic limit.}
\label{thermodynamic-limit}
\setcounter{equation}{0}

In this section, making use of the results of \cite{take:95}, we compute
several thermodynamic quantities of our spin chains.

First let us recall several facts on the solutions of the Bethe equations
\eqref{bethe-eq} that we found in \cite{take:95}. These solutions satisfy
the string hypothesis which goes back to Bethe \cite{bet}. For later
convenience we rescale the parameters as follows: 
$x = \i t u$, $x_j =\i t w_j$. In this notation the string hypothesis says
that for sufficiently large $N$ solutions of \eqref{bethe-eq} cluster into
groups known as $A$-strings, $A = 1, 2, \dots$, with parity $\pm$ and
centre $x^{A,\pm}_j$:
\begin{equation}
    x^{A,\pm}_{j,\alpha}= x^{A,\pm}_j + 2\i\eta t \alpha
                        + O(\e^{-\delta N}),
    \qquad
    \alpha = \frac{-A+1}{2}, \frac{-A+3}{2}, \dots, \frac{A-1}{2},
\label{def:A-string}
\end{equation}
where $\Im x^{A,+}_j = 0$ and $\Im x^{A,-}_j = t/2$. We denote the number
of $A$-strings with parity $\pm$ by $\sharp(A,\pm)$.

We consider the following string configurations which are consistent with
the constraints found in \cite{tak-suz} when $\eta = r'/r$, $r$, $r'$ are
integers mutually coprime, $r$ is even, $r'$ is odd, and 
$2(2l+1) \eta < 1$, as assumed in \cite{take:95}.

\begin{itemize}
\item Ground state:
$\nu=0$,
$\sharp(2l, +) = N/2$, $\sharp(A,\pm) = \sharp(2l, -) =0$
for $A\neq 2l$ and centres of $2l$-strings distribute
symmetrically around 0. 

\item Excited state I:
$\sharp(2l  , +) = N/2-2$,
$\sharp(2l-1, +) = 1$, 
$\sharp(2l+1, +) = 1$.

There are two holes in the distribution of $2l$-strings which are denoted
by $x_1$ and $x_2$ and regarded as continuous parameters of the
configuration. The Bethe equations determine the coordinates of the centres
of the $(2l\pm1)$-string which are denoted by $x_\pm$. There are two
possibilities: $(x_-, x_+) = ((x_1+x_2)/2, (x_1+x_2)/2)$ or
$(x_-, x_+) = ((x_1+x_2)/2, (x_1+x_2+1)/2)$.

\item Excited state II:
$\sharp(2l , +) = N/2-1$, 
$\sharp(2l-1, +) = 1$, 
$\sharp(1 , -) = 1$. 

There are again two holes in the distribution of $2l$-strings which are
denoted by $x_1$ and $x_2$. The Bethe equations determine the coordinates
of the centres of the $(2l-1)$-string and the 1-string with parity $-$
which are denoted by $x_-$ and $x_0$ respectively: 
$(x_-, x_0) = ((x_1+x_2)/2, (x_1+x_2)/2)$ or 
$(x_-, x_0) = ((x_1+x_2)/2, (x_1+x_2+1)/2)$.

\end{itemize}

Using the results of \cite{take:95}, we obtain the following asymptotics
of the eigenvalue of the transfer matrix $t^{l,l}(u)$ for large $N$. The
largest eigenvalue of $t^{l,l}$ which corresponds to the ground state is:
\begin{multline}
    \frac{1}{\i} \log t^{l,l}(u; \text{ ground state })
    \sim \\
    \sim
    N\left(
    \pi l + 2\pi l x (1-4l\eta) +
    2 \sum_{n=1}^\infty
    \frac{\sh 4\pi n l \eta t \sh \pi n t (1-4l\eta)}
         {n \sh \pi n t \sh 4 \pi n \eta t}
    \sin 2 \pi n x
    \right).
\label{transfer:gr-st}
\end{multline}
The excited states I and II have the same eigenvalue of the transfer
matrix $t^{l,l}(u)$. Namely,
\begin{multline}
    \frac{1}{\i} \log t^{l,l}(u;\, \text{excited state I or II})
    -
    \frac{1}{\i} \log t^{l,l}(u;\, \text{ground state})
    \\
    = \log \tau(x - x_1) + \log \tau(x - x_2),
\label{transfer:ex-st}
\end{multline}
where
\begin{equation}
    \log\tau(x) := - \frac{\pi}{2} - \pi x
                   - \sum_{n=1}^\infty
                   \frac{\sin 2\pi n x}{n \ch 2\pi n \eta t}.
\end{equation}
(see \cite{jo-kr-mc} and \cite{take:95} for details of computations). By
the definition \eqref{def:p,H}, the momentum and the energy of these
excited states are expressed as $p = p_1 + p_2$, $H = H_1 + H_2$ where
\begin{align}
    p_i &:= - \frac{\pi}{2} + \pi x_i
           + \sum_{n=1}^\infty
             \frac{\sin 2\pi n x_i}{n \ch 2\pi n \eta t},
\label{p:particle}
\\
    H_i &:= \text{const.}  \left(
            -\pi - 2 \pi \sum_{n=1}^\infty
            \frac{\cos 2\pi n x_i}{\ch 2\pi n \eta t}
            \right).
\label{H:particle}
\end{align}
Thus they are regarded as two particle spin wave modes each particle of
which has rapidity $x_i$, momentum $p_i$ and energy $H_i$. In particular
\eqref{p:particle} and \eqref{H:particle} show that the dispersion
relation of these particles do not depend on $l$.

\section*{Acknowledgements.}
The author expresses gratitude to Professor Nicolai Reshetikhin and
Professor Koji Hasegawa for discussions and comments and to the
Department of Mathematics of the University of California at Berkeley
for the hospitality. This work is partly supported by the Postdoctoral
Fellowship for Research abroad of Japan Society for the Promotion of
Science.

%
%
%
%

%
\end{document}